\documentstyle{amsppt}
\magnification=1200
\voffset-1cm
\NoRunningHeads
\def\0{\text{\bf 0}}
\def\cle{\preccurlyeq}
\def\cge{\succcurlyeq}
\def\Inf{$\wedge$}
\def\wo{$wo$}

\def\rset{\text{\bf R}}

\def\Rmax{{\text{\bf R}}_{\max}}
\def\Rmin{{\text{\bf R}}_{\min}}
\def\rmin{{\text{\bf R}}_{\min}}

\def\1{\text{\bf 1}}

\def\bR{\text{\bf R}}

\def\End{\text{\rm End}}

\define\pd#1#2{\dfrac{\partial#1}{\partial#2}}
\def\opl{\operatornamewithlimits{\oplus}}

\def\sumlim{\sum\limits}
\def\maxlim{\max\limits}

\def\suplim{\sup\limits}

\topmatter

\title{Idempotent (Asymptotic) Mathematics and the 
Representation Theory}\endtitle
\author {G.~L.~Litvinov, V.~P.~Maslov, and G.~B.~Shpiz}
\footnotemark"*"
\endauthor

\thanks
This work was supported by the RFBR grant 99--01--00196
and the Dutch Organization for Scientific Research (N.W.O.).
\endthanks

\endtopmatter

\document

\subhead 1. Introduction\endsubhead
Idempotent Mathematics is based on replacing the usual arithmetic
operations by a new set of basic operations (e.g., such as maximum
or minimum), that is on replacing numerical fields by idempotent 
semirings and semifields. Typical (and the most common) examples 
are given by the so-called $(\max, +)$ algebra $\Rmax$ and $(\min, +)$
algebra $\Rmin$. Let $\bR$ be the field of real numbers.
Then ${\Rmax}={\bR} \cup \{-\infty\}$ with operations
$x\oplus y=\max \{x,y\}$ and $x\odot y=x+y$. Similarly
${\Rmin}={\bR}\cup \{+\infty\}$ with the operations 
$\oplus=\min$, $\odot=+$. The new addition $\oplus$ is idempotent,
i.e., $x\oplus x=x$ for all elements $x$. Idempotent
Mathematics can be treated as a result of a dequantization of the
traditional mathematics over numerical fields as the Planck constant
$\hbar$ tends to zero taking pure imaginary values. Some problems
that are nonlinear in the traditional sense turn out to be linear
over a suitable idempotent semiring (idempotent superposition 
principle \cite{1}). For example, the Hamilton-Jacobi equation (which is an
idempotent version of the Schr\"odinger equation) is linear over
$\Rmin$ and $\Rmax$.

The basic paradigm is expressed in terms of an {\it idempotent 
correspondence principle} \cite{2}.This principle is similar to the 
well-known correspondence
principle of N.~Bohr in quantum theory (and closely related to it).
Actually, there exists a heuristic correspondence between
important, interesting and useful constructions and results of the
traditional mathematics over fields and analogous constructions and results
over idempotent semirings and semifields (i.e., semirings and semifields
with idempotent addition).

A systematic and consistent application of the idempotent
correspondence principle leads to a variety of results, often quite
unexpected. As a result, in parallel with the traditional mathematics over
rings, its ``shadow'', the Idempotent Mathematics, appears. This ``shadow''
stands approximately in the same relation to the traditional mathematics as
classical physics to quantum theory. In many respects Idempotent
Mathematics is simpler than the traditional one. However, the transition
from traditional concepts and results to their idempotent analogs is often
nontrivial.

There is an idempotent version of the theory of linear representations of 
groups. We shall present some basic concepts and results of the idempotent 
representation theory. In the framework of this 
theory the well-known Legendre transform can be treated as an 
$\bR_{\max}$-version of the traditional Fourier transform (this
observation is due to V.P. Maslov). We shall discuss some unexpected 
theorems of the  Engel type.

In this paper we present the so-called algebraic approach to Idempotent 
Mathematics: basic notions and results are `simulated' in pure
algebraic terms. Historical surveys and the corresponding references
can be found in \cite{2--6}.

\subhead{2. Semirings, semifields, and idempotent dequantization}\endsubhead
Consider a
set $S$ equipped with two algebraic operations: {\it addition} $\oplus$
and {\it multiplication} $\odot$. It is a {\it semiring} if the following
conditions are satisfied:
\roster
\item"{$\bullet$}" the addition $\oplus$ and the multiplication $\odot$ are
associative;
\item"{$\bullet$}" the addition $\oplus$ is commutative;
\item"{$\bullet$}" the multiplication $\odot$ is distributive with respect to 
the addition $\oplus$: $x\odot(y\oplus z)=(x\odot y)\oplus(x\odot z)$ and
$(x\oplus y)\odot z=(x\odot z)\oplus(y\odot z)$ for all $x,y,z\in S$.
\endroster
A {\it unity} of a semiring $S$ is an element $\1\in S$ such that
$\1\odot x=x\odot\1=x$ for all $x\in S$. A {\it zero} of a semiring
$S$ is an element $\0\in S$ such that $\0\neq\1$ and $\0\oplus x=x$, 
$\0\odot x=x\odot \0=\0$ for all $x\in S$. A semiring $S$ is called an
{\it idempotent semiring} if $x\oplus x=x$ for all $x\in S$. A 
semiring $S$ with neutral elements $\0$ and $\1$ is called a {\it 
semifield} if every nonzero element of $S$ is invertible.

Let $\rset$ be the field of real numbers and $\rset_+$ the semiring of all
nonnegative real numbers (with respect to the usual addition and 
multiplication).
The change of variables $x \mapsto u = h \ln x$, $h > 0$, defines a map
$\Phi_h \colon \rset_+ \to S = \rset \cup \{-\infty\}$.  Let the addition
and multiplication operations be mapped from $\rset$ to $S$ by $\Phi_h$,
i.e., let $u \oplus_h v = h \ln (\exp (u/h) + \exp(v/h))$, $u \odot v = u +
v$, $\0 = -\infty = \Phi_h(0)$, $\1 = 0 = \Phi_h(1)$.  It can easily be
checked that $u \oplus_h v \to \max \{u,v\}$ as $h \to 0$ and $S$ forms a
semiring with respect to addition $u \oplus v = \max \{u,v\}$ and
multiplication $u \odot v = u + v$ with zero $\0 = -\infty$ and unit $\1 =
0$. Denote this semiring by $\Rmax$; it is {\it idempotent}, i.e., $u
\oplus u = u$ for all its elements.
The semiring $\bR_{\max}$ is actually a semifield.
The analogy with quantization is obvious; the parameter $h$ plays the
r\^{o}le of the Planck constant, so $\rset_+$ (or $\rset$) can be viewed as
a ``quantum object'' and $\Rmax$ as the result of its ``dequantization''. A
similar procedure gives the semiring $\rmin = \rset \cup \{+\infty\}$ with
the operations $\oplus = \min$, $\odot = +$; in this case $\0 = +\infty$, $\1 =
0$. The semirings $\Rmax$ and $\rset_{\min}$ are isomorphic. Connections
with physics and imaginary values of the Planck constant are discussed
below. The idempotent semiring $\rset \cup \{-\infty\} \cup
\{+\infty\}$ with the operations $\oplus = \max$, $\odot = \min$ can be
obtained as a result of a ``second dequantization'' of $\rset$ (or
$\rset_+$).  Dozens of interesting examples of nonisomorphic idempotent
semirings may be cited as well as a number of standard methods of deriving
new semirings from these (see, e.g., \cite{2--6} and below).

Idempotent dequantization is related to logarithmic transformations that
were used in the classical papers of E. Schr\"odinger \cite{7} and
E. Hopf \cite{8}. The subsequent progress of E. Hopf's ideas has
culminated in the well-known vanishing viscosity method (the method
of viscosity solutions), see, e.g., \cite{9}.
 
\subhead 3. Idempotent Analysis\endsubhead
Let $S$ be an arbitrary semiring with idempotent addition $\oplus$ (which
is always assumed to be commutative), multiplication $\odot$, zero $\0$, and
unit $\1$. The set $S$ is supplied with the {\it standard partial
order\/}~$\cle$: by definition, $a \cle b$ if and only if $a \oplus b = b$.
Thus all elements of $S$ are positive: $\0 \cle a$ for all $a \in S$. Due
to the existence of this order, Idempotent Analysis is closely related to
lattice theory, the theory of vector lattices, and the theory
of ordered spaces. Moreover, this partial order allows to
model a number of basic notions and results of Idempotent Analysis at
the purely algebraic level; in this paper we develop
this line of reasoning systematically. Let us notice that the standard
partial order can be defined for an arbitrary commutative semigroup with
idempotent addition.

Calculus deals mainly with functions whose values are numbers. The
idempotent analog of a numerical function is a map $X \to S$, where $X$ is
an arbitrary set and $S$ is an idempotent semiring. Functions with values
in $S$ can be added, multiplied by each other, and multiplied by elements
of $S$.

The idempotent analog of a linear functional space is a set of $S$-valued
functions that is closed under addition of functions and multiplication of
functions by elements of $S$, or an $S$-semimodule. Consider, e.g., the
$S$-semimodule $\Cal B(X, S)$ of functions $X \to S$ that are bounded in
the sense of the standard order on $S$.

If $S = \Rmax$, then the idempotent analog of integration is defined by the
formula
$$
   I(\varphi) = \int_X^{\oplus} \varphi (x)\, dx
        = \sup_{x\in X} \varphi (x),\tag{1}
$$
where $\varphi \in \Cal B(X, S)$. Indeed, a Riemann sum of the form
$\sumlim_i \varphi(x_i) \cdot \sigma_i$ corresponds to the expression
$\bigoplus\limits_i \varphi(x_i) \odot \sigma_i = \maxlim_i \{\varphi(x_i)
+ \sigma_i\}$, which tends to the right-hand side of~(1) as $\sigma_i \to
0$. Of course, this is a purely heuristic argument.

Formula~(1) defines the idempotent integral not only for functions taking
values in $\Rmax$, but also in the general case when any of bounded
(from above) subsets of~$S$ has the least upper bound.

An idempotent measure on $X$ is defined by $m_{\psi}(Y) = \suplim_{x \in Y}
\psi(x)$, where $\psi \in \Cal B(X,S)$. The integral with respect to this
measure is defined by
$$
   I_{\psi}(\varphi)
        = \int^{\oplus}_X \varphi(x)\, dm_{\psi}
        = \int_X^{\oplus} \varphi(x) \odot \psi(x)\, dx
        = \sup_{x\in X} (\varphi (x) \odot \psi(x)).
        \tag{2}
$$

Obviously, if $S = \rmin$, then the standard order $\cle$ is opposite to
the conventional order $\leqslant$, so in this case equation~(2) assumes the form
$$
   \int^{\oplus}_X \varphi(x)\, dm_{\psi}
        = \int_X^{\oplus} \varphi(x) \odot \psi(x)\, dx
        = \inf_{x\in X} (\varphi (x) \odot \psi(x)),
        \tag{3}
$$
where $\inf$ is understood in the sense of the conventional order $\leqslant$.

The functionals $I(\varphi)$ and $I_{\psi}(\varphi)$ are linear over $S$;
their values correspond to limits of Lebesgue (or Riemann) sums.
The formula for $I_\psi(\varphi)$ defines the
idempotent scalar product of the functions $\psi$ and $\varphi$. Various
idempotent functional spaces and an idempotent version of the theory of
distributions can be constructed on the basis of the idempotent integration,
see, e.g., \cite{1, 3--6, 10}. The analogy of idempotent and probability measures leads to
spectacular parallels between optimization theory and probability theory. For
example, the Chapman--Kolmogorov equation corresponds to the Bellman
equation (see, e.g., the survey of Del~Moral \cite{11} and \cite{5}).
Many other idempotent analogs may be cited (in particular, for basic
constructions and theorems of functional analysis \cite{4}).

\subhead 4. The superposition principle and linear problems\endsubhead
 Basic equations of quantum theory
are linear (the superposition principle). The Hamilton--Jacobi equation, the basic equation of classical
mechanics, is nonlinear in the conventional sense. However it is linear over the
semiring $\rmin$. Also, different versions of the Bellman equation, the basic
equation of optimization theory, are linear over suitable idempotent
semirings (V.~P.~Maslov's idempotent superposition
principle), see \cite{1, 3, 6, 10}. For instance, the finite-dimensional
stationary Bellman equation can be written in the form $X = H \odot
X \oplus F$, where $X$, $H$, $F$ are matrices with coefficients in
an idempotent semiring $S$ and the unknown matrix $X$ is determined by $H$
and $F$ \cite{12}. In particular, standard problems of dynamic programming
and the well-known shortest path problem correspond to the cases $S = \Rmax$ 
and $S =\rmin$, respectively. In \cite{12}, it was shown that main 
optimization algorithms for finite graphs correspond to standard methods for 
solving systems of linear equations of this type (i.e., over semirings).
Specifically, Bellman's shortest path algorithm
corresponds to a version of Jacobi's algorithm, Ford's algorithm
corresponds to the Gauss--Seidel iterative scheme, etc.

Linearity of the Hamilton--Jacobi equation over $\rmin$ (and $\Rmax$) is
closely related to the (conventional) linearity of the Schr\"odinger equation.
Consider a classical dynamical system specified by the Hamiltonian
$$
   H = H(p,x) = \sum_{i=1}^N \frac{p^2_i}{2m_i} + V(x),\tag{6}
$$
where $x = (x_1, \dots, x_N)$ are generalized coordinates, $p = (p_1,
\dots, p_N)$ are generalized momenta, $m_i$ are generalized masses, and
$V(x)$ is the potential. In this case the Lagrangian $L(x, \dot x, t)$ has
the form
$$
   L(x, \dot x, t)
        = \sum^N_{i=1} m_i \frac{\dot x_i^2}2 - V(x),\tag{7}
$$
where $\dot x = (\dot x_1, \dots, \dot x_N)$, $\dot x_i = dx_i / dt$. The
value function $S(x,t)$ of the action functional has the form
$$
   S(x, t)= \int^t_{t_0} L(x(t), \dot x(t), t)\, dt,\tag{8}
$$
where the integration is performed along a trajectory of the system.  The
classical equations of motion are derived as the stationarity conditions
for the action functional (the Hamilton principle, or the least action
principle).

The action functional can be considered as a function taking the set of
curves (trajectories) to the set of real numbers. Assume that its range
lies in the semiring $\rmin$. In this case the minimum of the action
functional can be viewed as the idempotent integral of this function over
the set of trajectories or the idempotent analog of the Feynman path
integral. Thus the least action principle can be considered as the idempotent
version of the well-known Feynman approach to quantum mechanics (which is
presented, e.g., in~\cite{13}); here, one should remember that the
exponential function involved in the Feynman integral is monotone on the
real axis. The representation of a solution to the Schr\"{o}dinger
equation in terms of the Feynman integral corresponds to the
Lax--Ole\u{\i}nik formula for a solution to the Hamilton--Jacobi equation.

Since $\partial S/\partial x_i = p_i$, $\partial S/\partial t = -H(p,x)$,
the following Hamilton--Jacobi equation holds:
$$
   \pd{S}{t} + H \left(\pd{S}{x_i}, x_i\right)= 0.\tag{9}
$$

Quantization (see, e.g., \cite{13}) leads to the Schr\"odinger equation
$$
   -\frac{\hbar}i \pd{\psi}{t}= \widehat H \psi = H(\hat p_i, \hat x_i)\psi,
   \tag{10}
$$
where $\psi = \psi(x,t)$ is the wave function, i.e., a time-dependent
element of the Hilbert space $L^2(\rset^N)$, and $\widehat H$ is the energy
operator obtained by substitution of the momentum operators $\widehat p_i
= {\hbar \over i}{\partial \over \partial x_i}$ and the coordinate
operators $\widehat x_i \colon \psi \mapsto x_i\psi$ for the variables
$p_i$ and $x_i$ in the Hamiltonian function, respectively. This equation is
linear in the conventional sense (the quantum superposition principle). The
standard procedure of limit transition from the Schr\"odinger equation to
the Hamilton--Jacobi equation is to use the following ansatz for the wave
function:  $\psi(x,t) = a(x,t) e^{iS(x,t)/\hbar}$, and to keep only the
leading order as $\hbar \to 0$ (the `semiclassical' limit).

Instead of doing this, we switch to imaginary values of the Planck constant
$\hbar$ by the substitution $h = i\hbar$, assuming $h > 0$. Thus the
Schr\"odinger equation~(1.10) turns to an analog of the heat equation:
$$
   h\pd{u}{t} = H\left(-h\frac{\partial}{\partial x_i}, \hat x_i\right) u,
   \tag{11}
$$
where the real-valued function $u$ corresponds to the wave function
$\psi$. A similar idea (the switch to imaginary time) is used in the
Euclidean quantum field theory (see, e.g., \cite{14}); let us remember
that time and energy are dual quantities.

Linearity of equation~(10) implies linearity of equation~(11). Thus if
$u_1$ and $u_2$ are solutions of~(11), then so is their linear
combination
$$
   u = \lambda_1 u_1 + \lambda_2 u_2.\tag{12}
$$

Let $S = -h \ln u$ or $u = e^{-S/h}$ as in Section 2 above. It can easily
be checked that equation~(11) thus turns to
$$
   \pd{S}{t}= V(x) + \sum^N_{i=1} \frac1{2m_i}\left(\pd{S}{x_i}\right)^2
        - h\sum^n_{i=1}\frac1{2m_i}\frac{\partial^2 S}{\partial x^2_i}.
   \tag{13}
$$
This equation is nonlinear in the conventional sense. However, if $S_1$ and $S_2$
are its solutions, then so is the function
$$
   S = \lambda_1 \odot S_1 \opl_h \lambda_2\odot S_2\tag{14}
$$
obtained from~(12) by means of our substitution $S = -h \ln u$.  Here the
generalized multiplication $\odot$ coincides with the ordinary addition and
the generalized addition $\opl_h$ is the image of the conventional addition
under the above change of variables.  As $h \to 0$, we obtain the
operations of the idempotent semiring $\rmin$, i.e., $\oplus = \min$ and
$\odot = +$, and equation~(13) turns to the Hamilton--Jacobi
equation~(9), since the third term in the right-hand side of
equation~(13) vanishes.

Thus it is natural to consider the limit function $S = \lambda_1 \odot S_1
\oplus \lambda_2 \odot S_2$ as a solution of the Hamilton--Jacobi equation and
to expect that this equation can be treated as linear over $\rmin$. This
argument (clearly, a heuristic one) can be extended to equations of a more
general form. For a rigorous treatment of (semiring) linearity for these
equations see \cite{3, 6} and also \cite{1}. Notice that if $h$ is changed
to $-h$, then the resulting Hamilton--Jacobi equation is linear over
$\Rmax$.

The idempotent superposition principle indicates that there exist important
problems that are linear over idempotent semirings.

\subhead{5. Convolution and the Fourier--Legendre transform}\endsubhead 
Let $G$ be a group. Then the space $\Cal B(X, \bR_{\max})$ of all bounded
functions $G\to\bR_{\max}$ (see above) is an idempotent semiring
with respect to the following analog $\circledast$ of the usual 
convolution:
$$
   (\varphi(x)\circledast\psi)(g)=
        = \int_G^{\oplus} \varphi (x)\odot\psi(x^{-1}\cdot g)\, dx=
\sup_{x\in G}(\varphi(x)+\psi(x^{-1}\cdot g)).\tag{15}
$$
Of course, it is possible to consider other ``function spaces'' (and other
basic semirings instead of $\bR_{\max}$). In \cite{3} ``group
semirings'' of this type are referred to as {\it convolution
semirings}.

Let $G=\bR^n$, where $\bR^n$ is considered as a topological group
with respect to the vector addition. The conventional Fourier--Laplace 
transform is defined as
$$
   \varphi(x) \mapsto \tilde{\varphi}(\xi)
        = \int_G e^{i\xi \cdot x} \varphi (x)\, dx,\tag{16}
$$
where $e^{i\xi \cdot x}$ is a character of the group $G$, i.e., a solution
of the following functional equation:
$$
   f(x + y) = f(x)f(y).
$$
The idempotent analog of this equation is
$$
   f(x + y) = f(x) \odot f(y) = f(x) + f(y),
$$
so ``continuous idempotent characters'' are linear functionals of the
form $x \mapsto \xi \cdot x = \xi_1 x_1 + \dots + \xi_n x_n$. As a result,
the transform in~(16) assumes the form
$$
   \varphi(x) \mapsto \tilde{\varphi}(\xi)
        = \int_G^\oplus \xi \cdot x \odot \varphi (x)\, dx
   = \sup_{x\in G} (\xi \cdot x + \varphi (x)).\tag{17}
$$
The transform in~(17) is nothing but the {\it Legendre transform\/}
(up to some notation) \cite{10}; transforms of this kind establish the
correspondence between the Lagrangian and the Hamiltonian formulations
of classical mechanics.

Of course, this construction can be generalized to different classes
of groups and semirings. Transformations of this type convert the
generalized convolution $\circledast$ to the pointwise (generalized)
multiplication and possess analogs of some important properties
of the usual Fourier transform. For the case of semirings of Pareto
sets the corresponding version of the Fourier transform reduces the
multicriterial optimization problem to a family of singlecriterial 
problems \cite{15}.

The examples discussed in this sections can be treated as fragments
of an idempotent version of the representation theory. In particular,
``idempotent'' representations of groups can be examined as representations
of the corresponding convolution semirings (i.e. idempotent group semirings)
in semimodules. To present nontrivial examples from the idempotent 
version of the representation theory we need some preliminary material.

\subhead{6. Idempotent semimodules and linear spaces}\endsubhead
Recall that an idempotent semigroup is an arbitrary commutative 
(additive) semigroup with idempotent addition. It can be treated
as an ordered set with the following partial order: $x\cle y$ if
and only if $x\oplus y=y$. It is easy to see that 
this order is well-defined and $x\oplus y=\sup \{x,y\}$. For an arbitrary
subset $X$ of an idempotent semigroup, we put $\oplus X=\sup (X)$ 
and $\wedge X=\inf (X)$ if the corresponding right-hand sides exist.
An idempotent semigroup is called $b$-{\it complete} (or {\it boundedly 
complete}) if  any of its subsets
bounded from above (including the empty subset) has the least upper
bound. In particular, any $b$-complete idempotent semigroup contains zero 
(denoted by $\0 $), which coincides with $\oplus \varnothing$, where 
$\varnothing$ is the empty set. A homomorphism of $b$-complete idempotent 
semigroups is called a $b$-{\it homomorphism} if $g(\oplus X)=\oplus g(X)$ 
for any subset $X$ bounded from above.

An idempotent semifield is called $b$-{\it complete} if it is
$b$-complete as an idempotent semigroup.
In any $b$-complete semifield, the generalized distributive laws
$$
a\odot (\oplus X)=\oplus (a\odot X),\quad 
a\odot (\wedge X)=\wedge (a\odot X)\tag{18}
$$
are valid; here $a$ is an element of the semifield and $X$ is a nonempty 
bounded subset. It is easy to see that $\bR_{\max}$ is a $b$-complete
semifield.

An {\it idempotent semimodule} over an idempotent semiring $K$ is an 
idempotent semigroup $V$ endowed with a multiplication $\odot$ by elements 
of $K$ such that, for any $a,b\in K$ and $x,y\in V$, the usual laws
$$
\align
a\odot (b\odot x)&=(a\odot b)\odot x, \tag{19}\\
(a\oplus b)\odot x&=a\odot x\oplus b\odot x,\tag{20}\\
a\odot (x\oplus y)&=a\odot x\oplus a\odot y, \tag{21} \\ 
\0 \odot x&=\0\tag{22}
\endalign
$$
are valid. An idempotent semimodule over an idempotent semifield is called
an {\it idempotent space}. An idempotent $b$-complete space $V$ 
over a $b$-complete semifield
$K$ is called an {\it idempotent} $b$-{\it space} if, for
 any nonempty bounded subset 
$Q\subset K$ and any $x\in V$, the relations
$$
(\oplus Q)\odot x=\oplus (Q\odot x),\quad 
(\wedge Q)\odot x=\wedge (Q\odot x)\tag{23}
$$
hold. A homomorphism $g: V\to W$ of $b$-spaces is called a $b$-{\it
homomorphism}, or a $b$-{\it linear operator} ({\it mapping}), 
if $g(\oplus X)=\oplus g(X)$ 
for any bounded subset $X\subset V$. More general definitions 
(for spaces which may not be $b$-complete) can be found in \cite{4}.
Homomorphisms taking values in $K$ (treated as a
semimodule over itself) are called {\it linear functionals}. A subset of an
idempotent space is called a {\it subspace} if it is closed with respect to
addition and multiplication by coefficients. A subspace in a $b$-space is 
called a $b$-{\it closed subspace} if it is closed with respect to 
summation over arbitrary bounded (in V) subsets. This subspace has a natural
structure of $b$-space; it is also a $b$-subspace in V in the sense of 
\cite{4}.

For an arbitrary set $X$ and an idempotent space $V$ over a semifield $K$,
we use $B(X,V)$ to denote the semimodule of all bounded mappings from $X$ 
into $V$ with pointwise operations. If $V$ is an idempotent $b$-space,
then $B(X,V)$ is a $b$-space. A mapping $f$ from a topological space $X$ 
into an ordered set $V$ is called {\it upper semicontinuous} if, for any
$b\in V$, the set  $\{x\in X | f(x)\cge b\}$ is closed in $X$, see \cite{4}. 
In the case where $V$ is the set of real numbers, this definition 
coincides with the usual
definition of upper semicontinuity of a real function. The set of all
bounded upper semicontinuous mappings from $X$ to $V$ is denoted by 
$USC(X,V)$. If V is an idempotent $b$-space, then $USC(X,V)$ is also a 
$b$-space with respect to the operations $f\oplus g=\sup \{f,g\}$ and 
$(k\odot f)(x)=k\odot f(x)$. 

\subhead{7. Archimedean spaces \cite{16}}\endsubhead
In what follows, unless otherwise specified, the symbol $K$
stands for a $b$-complete idempotent semifield and all idempotent spaces
are over  $K$. 

A subset $M$ of idempotent $b$-space $V$ is called \wo-{\it closed} if 
$\wedge X\in M$ and $\oplus X\in M $ for any linearly ordered subset 
$X\nomathbreak\subset\nomathbreak M$ in $V$. A nondecreasing mapping
$f: V\to W$ of $b$-spaces is called \wo-continuous if  
$f(\oplus X)=\oplus f(X)$ and $f(\wedge X)= \wedge f(X)$ for any bounded
linearly ordered subset $X\subset V$. Note that an arbitrary isomorphism of
ordered sets is \wo-continuous. It can be shown that the notions of 
\wo-closedness and \wo-continuity coincide with the closedness and continuity 
with respect to some $T1$ topology defined in an intrinsic way in terms of the
order.

\proclaim{Proposition 1} Suppose that $V$ is an idempotent b-space and 
$W$ is a \wo-closed subsemigroup of $V$. Then $\oplus \nomathbreak X
\in \nomathbreak W$ for any subset $X\subset W$ bounded in $V$. In
particular, each \wo-closed subspace is a b-closed subspace.\endproclaim

An element $x$ of an idempotent space $V$ is called {\it Archimedean} if, 
for any $y\in V$, there exists a coefficient $\lambda \in K$ such that  
$\lambda \odot x\cge y$. For an Archimedean element $x\in V$, the
formula $x^*(y)=\wedge \{k\in K | k\odot x\cge y\}$ defines a mapping 
$x^*: V\to K$. If $V$ is an idempotent $b$-space, then $x^*$ is a
$b$-linear functional and $x^*(y)\odot x\cge y$ for any $y\in V$ [6]. 
We say that an Archimedean element $x\in\nomathbreak V$ is \wo-{\it continuous}
if the functional $x^* $ is \wo-continuous, and that an idempotent 
$b$-space $V$ is {\it Archimedean} if $V$ contains a \wo-continuous Archimedean
element.

\proclaim{Proposition 2} If $X$ is a compact topological space, then 
USC(X,K) is an Archi\-medean space and the function $\bold e$ identically
equal to $\1$ is a \wo-continuous Archime\-dean element.\endproclaim

Note that $\bold {e^*}(f)=\sup \{f(x) | x\in X\}$.

\proclaim{Theorem 1} Any \wo-closed subspace of an Archimedean space
is an Archimedean space. Any linearly ordered (with respect to the
inclusion) family of nonzero \wo-clos\-ed subspaces of an Archimedean
space $V$ has a nonzero intersection.\endproclaim

 Let $V$ be a $b$-space. A subset $W\subset V$ is called a 
\Inf-{\it subspace} if it is closed with respect to multiplication by scalars
and taking greatest lower bounds of nonempty subsets. By this definition, any
such $W$ is a boundedly complete lattice with respect to the order inherited
from $V$. Therefore, any \Inf-subspace $W\subset V$ can
be treated as a semimodule with respect to the inherited multiplication
by scalars and the operations $x\oplus _W y=\sup \{x,y\}$, where $\sup$ 
is over $W$. In what  follows, all \Inf-subspaces are considered as 
semimodules with respect to these operations. The definitions
immediately imply that any \Inf-subspace of a $b$-space is a $b$-space. 
It is easy to show that $USC(X,V)$ is a \Inf-subspace in $B(X,V)$ 
for any $b$-space $V$ and any topological space $X$. 

\proclaim{Proposition 3} If $V$ is an Archimedean $b$-space and $x\in V$ 
is a \wo-continuous Archi\-medean element, then any \Inf-subspace $W$ 
of $V$ containing $x$ is an Archimedean $b$-space.\endproclaim

An arbitrary semiring $K$ is called {\it algebraically closed}
(or {\it radicable}, see, e.g. \cite{5}) if for any element
$x\in K$ and any positive integer number $n$ there exists an
element $y\in K$ such that $y^n=x$. It is easy to show that $\bR_{\max}$
is a $b$-complete algebraically closed semifield.

\proclaim{Theorem 2} An idempotent $b$-space $V$ over an algebraically
closed $b$-complete semifield $K$ is Archimedean if and only if there exists 
a space of the form USC{\rm(}X,K{\rm)}, where $X$ is a compact topological space,
such that $V$ is isomorphic to its \Inf-subspace containing constants.
\endproclaim

\subhead{8. Representations of groups in Archimedean spaces}\endsubhead
Suppose that $V$ is an Archimedean idempotent $b$-space over an 
algebraically closed $b$-complete semifield (e.g., over $\bR_{\max}$).
By $\End (V)$ denote the set of all $b$-linear operators $V\to V$.
This set is an idempotent semigroup with respect to the pointwise
sum and it is a $b$-space over $K$ with respect to the standard
multiplication by coefficients from $K$. The usual multiplication
(composition) of maps turns $\End (V)$ into an idempotent semiring
(and a $b$-complete semialgebra over $K$).

Let $G$ be an abstract group. A {\it linear representation}
$\pi : G \to \End (V)$ of $G$ in an Archimedean space $V$ is a homomorphism
from $G$ to the group of all invertible elements in $\End (V)$ (with
respect to the composition of operators). The representation $\pi$
is (topologically) {\it irreducible} if the space $V$ has no nontrivial
$wo$-closed $\pi (G)$-invariant subspaces.

\proclaim{Theorem 3} Every linear representation of a group $G$
in an Archimedean idempotent space $V$ has a nontrivial irreducible
subrepresentation in a $wo$-closed subspace of $V$.\endproclaim

\proclaim{Theorem 4} Let $\pi$ be a linear representation of a group $G$
in an Archimedean idempotent space $V$ and for a nonzero element
$x\in V$ the orbit $\pi(G)x$ is bounded. Set $a=\oplus(\pi(G)x)$.
Then $\pi(g)a=a$ for each $g\in G$.\endproclaim

We shall say that a representation $\pi$ of $G$ in $V$ has a (nonzero)
{\it joint eigenvector} $a\in V$ if $\pi(g)a=\lambda(g)a$ for all
$g\in G$, where $\lambda(g)\in K$.

\proclaim{Corollary 1} Every linear representation of a finite group
in an Archimedean idempotent space has a joint eigenvector with a
unique eigenvalue $\1$.\endproclaim

\proclaim{Corollary 2} Every upper semicontinuous linear representation 
of a compact group in an Archimedean idempotent space has a joint
eigenvector with a unique eigenvalue $\1$.\endproclaim

\subhead{9. An Engel type theorem for representations of nilpotent groups}
\endsubhead
Let $G$ be an abstract group. For elements $a, b \in G$ we set
$[a, b] = a^{-1}b^{-1} ab$; for subsets $X$ and $Y$ in $G$ we denote by 
$[X,Y]$ a subgroup in $G$ generated by the set $\{[x, y]  \vert x\in X,
y\in Y\}$; we set $\Gamma_1(G)=G$, $\Gamma_i(G)=[G, \Gamma_{i-1}(G)]$,
$i=1, 2, 3, \dots$ Recall that an abstract group $G$ is {\it nilpotent}
if and only if there exists a positive integer number $n$ such that
$\Gamma_n(G)=\{e\}$, where $e$ is the neutral element (identity) of $G$.

\proclaim{Theorem 5} Every linear representation of a nilpotent abstract
group in an Archi\-medean idempotent space over an algebraically closed
$b$-complete semifield (e.g., over $\bR_{\max}$) has a joint eigenvector.
\endproclaim

\proclaim{Corollary 3} Every collection of commuting invertible $b$-linear
operators in an Archimedean idempotent space has a joint eigenvector.
\endproclaim

\proclaim{Corollary 4} Every invertible $b$-linear operator in an 
arbitrary Archimedean idempotent space over an algebraically closed 
$b$-complete semifield has an eigenvector.\endproclaim

{\bf Remark.} There is no idempotent version of the well known Lie
theorem for representations of abstract solvable groups in 
idempotent spaces. Moreover, there exists an irreducible linear 
representation of a solvable group in the idempotent space
$V=\bR_{\max}\times\bR_{\max}$ over $\bR_{\max}$.

\enddocument